\documentclass[a4paper,12pt]{amsart}

\usepackage{amsthm, amsmath, latexsym, amsfonts,hyperref}

\theoremstyle{plain}
\newtheorem{Lemma}{Lemma}
\newtheorem{Proposition}{Proposition}

\newtheorem{Corollary}{Corollary}
\theoremstyle{definition}
\newtheorem{Remark}{Remark}

\newcommand{\monbar}{\overline{M}_{0,n}}
\newcommand{\monp}{\overline{M}_{0,n}(\mathbb{P}^1,1)}

\title{When is $\monp$ a Mori dream space?}

\author{Claudio Fontanari}

\email{claudio.fontanari@unitn.it}\curraddr{{\sc Dipartimento di Matematica \\  Universit\`a degli Studi di Trento\\ Via Sommarive 14 \\ 38123 Trento \\ Italy.}}

\thanks{This research was partially supported by PRIN 2017 ``Moduli Theory and Birational Classification" and by GNSAGA of INdAM (Italy).}



\begin{document}

\begin{abstract}
We prove that the moduli space of $n$-pointed stable maps $\monp$ is a Mori dream space whenever the moduli space $\overline{M}_{0,n+3}$ 
of $(n+3)$-pointed rational curves is. We also show that $\monp$ is a log Fano variety for $n \le 5$.
\end{abstract}

\maketitle

\section{Introduction}
Given a normal projective $\mathbb{Q}$-factorial variety $X$ over an algebraically closed field $\mathbb{K}$ of any
characteristic, $X$ is a Mori dream space if $X$ has the following properties (see for instance \cite{C:18}): (1) The Picard group 
$\mathrm{Pic}(X)$ of $X$ is finitely generated, and $\mathrm{Pic}(X)_\mathbb{Q} = \mathrm{N_1}(X)_\mathbb{Q}$; (2) the cone 
of nef divisors $\mathrm{Nef}(X)$ is generated by a finite number of semi-ample divisors; and (3) there are finitely many small,
$\mathbb{Q}$-factorial modifications $f_i : X  \dashrightarrow X_i$ of $X$ such that each $X_i$ has properties (1) and (2) and the
moving cone of $X$ is the union of the pullbacks of the nef cones of the $X_i$. In other words, if $X$ is a Mori dream space, 
then one would be able, at least in principle, to explicitly describe the birational models $X_i$ of $X$, which are isomorphic to $X$ in codimension one, and to use them to describe the nef cone of $X$. It would also follow that the effective cone of $X$ is polyhedral. 
The nef and effective cones of divisors are crucial in understanding the birational geometry of a variety. In particular, for moduli spaces of curves an understanding of how these cones relate to each other was a crucial ingredient in the proof that the moduli space of stable curves $\overline{M}_g$ is of general type for $g = 22$ and $g \ge 24$ (see for instance \cite{F:09}). Even for $g=0$ the moduli space 
$\monbar$, parameterizing stable rational curves with $n$ ordered marked points and not too far from being a toric variety, presents a surprisingly rich birational geometry. The partial results obtained in two decades of intensive investigation range from the positive side 
(for instance, $\monbar$ is a Mori dream space for $n \le 6$, see \cite{KM:96} and \cite{C:09}) to the negative one (as the breakthrough in \cite{CT:15} that $\monbar$ is not a Mori dream space for $n \ge 134$, later improved to $n \ge 13$ in \cite{GK:16} and then to $n \ge 10$ in \cite{HKL:18}). 

Here we address the same question for a different but closely related moduli space. As it is well-known (see for instance \cite{Ma:17}, Remark 1.4), the Kontsevich 
moduli space $\monp$ parameterizing $n$-pointed stable maps to $\mathbb{P}^1$ of genus $0$ and degree $1$ is isomorphic to the Fulton-MacPherson compactification 
$\mathbb{P}^1[n]$ of the configuration space of $n$ distinct ordered points in $\mathbb{P}^1$. The natural projection $\monp \to \monbar$ implies by \cite{O:16}
that if $\monp$ is a Mori dream space then $\monbar$ is a Mori dream space too. In particular, it follows that $\monp$ is not a Mori dream space for $n \ge 10$. 

In Section \ref{birational} we establish a converse statement: if $\overline{M}_{0,n+3}$ is a Mori dream space then $\monp$ is a Mori dream space too 
(see Proposition \ref{birmap}). In order to do so, we introduce a natural birational map $\overline{M}_{0,n+3} \dashrightarrow \monp$ which is surjective in codimension one 
and we apply \cite{O:15}. 
In particular, from the known results for $\monbar$ we recover the fact that $\monp$ is a Mori dream space for $n \le 3$, which is already well understood: indeed, $\mathbb{P}^1[1] \cong \mathbb{P}^1$, $\mathbb{P}^1[2] \cong \mathbb{P}^1 \times \mathbb{P}^1$ and 
$\mathbb{P}^1[3]$ appears in the list of smooth Fano threefolds  (for instance \cite{Ma:17}, p. 108), so it is a Mori dream space by \cite{BCHM:10}, Corollary 1.3.2. 

To go further we need to implement a different strategy. After reph\-rasing in Section~\ref{ample} the characterization of ample divisors on $\monp$ provided by 
\cite{CHS:09}, in Section~\ref{fano} we check that $\monp$ is a log Fano variety for $n \le 5$ but not for $n=6$. We conclude that $\monp$ is a 
Mori dream space for $n \le 5$ (see Corollary \ref{dream}) and we point out that new ideas are required to address the remaining open cases 
$6 \le n \le 9$ (see Remark \ref{nohope}). 

We work over the complex field $\mathbb{C}$.

We are grateful to the anonymous referee for detailed suggestions in order to improve the above Introduction. 

\section{}\label{birational}
First we recall the definition and the basic properties of both $\monbar$ and $\monp$ following \cite{FP:97}.

The moduli space $\monbar$ parameterizes isomorphism classes of stable curves of genus $0$ with $n$ ordered marked points: 
$$
(C, p_1, \ldots, p_n). 
$$
For every subset $S \subset \{ 1, \ldots, n \}$ with $2 \le \vert S \vert \le n-2$ the boundary component $\Delta_S$ is the closure in $\monbar$ of the locus of stable curves 
$$
(C_1= \mathbb{P}^1, (p_i)_{i \in S}) \cup (C_2= \mathbb{P}^1, (p_i)_{i \in S^c}).
$$ 
The moduli space $\monp$ parameterizes isomorphism classes of stable maps of degree $1$ from curves of genus $0$ with $n$ ordered marked points to $\mathbb{P}^1$: 
$$
(C, p_1, \ldots, p_n, f: C \to \mathbb{P}^1).
$$
For every subset $S \subset \{ 1, \ldots, n \}$ with $2 \le \vert S \vert \le n$ the boundary component $B_S$ is the closure in $\monp$ of the locus of stable maps:
$$
(C_1= \mathbb{P}^1, (p_i)_{i \in S}, \mathrm{pt}: \mathbb{P}^1 \to \mathbb{P}^1) \cup (C_2= \mathbb{P}^1, (p_i)_{i \in S^c}, \mathrm{id}: \mathbb{P}^1 \to \mathbb{P}^1)
$$
collapsing the first component to the point $C_1 \cap C_2$ and mapping the second component identically to $\mathbb{P}^1$.

Both $\monbar$ and $\monp$ are smooth projective varieties and in both cases the union of the boundary components is a normal crossing divisors (see for instance \cite{HT:03}, Theorem 2.3).  

\begin{Proposition}\label{birmap}  If $\overline{M}_{0,n+3}$ is a Mori dream space then $\monp$ is a Mori dream space.
\end{Proposition}

\proof 
By \cite{O:15}, Proposition 1.3 and Remark 2.2, the claim follows if there is a birational map
$\overline{M}_{0,n+3} \dashrightarrow \monp$ 
which is surjective in codimension one. 

Let 
$$
U_0 := \{ (C_0 \cup \ldots \cup C_k, p_1, \ldots, p_{n+3}) \in \overline{M}_{0,n+3}: p_{n+1}, p_{n+2}, p_{n+3} \in C_0 \}
$$
and notice that $U_0$ is dense in $\overline{M}_{0,n+3}$ since it contains the open part $M_{0,n+3} \subset \overline{M}_{0,n+3}$ 
parameterizing smooth rational curves. 

Consider the natural rational map: 
\begin{eqnarray*}
\Phi: U_0 \dashrightarrow \monp \\
(C_0= \mathbb{P}^1, (p_i)_{i \in S_0}, p_{n+1}, p_{n+2}, p_{n+3}) \cup
\bigcup_{j=1}^k (C_j= \mathbb{P}^1, (p_i)_{i \in S_j}) \\
\mapsto (C_0= \mathbb{P}^1, (\pi(p_i))_{i \in S_0}, \mathrm{id}) \cup 
\bigcup_{j=1}^k (C_j= \mathbb{P}^1, (p_i)_{i \in S_j}, \mathrm{pt}) 
\end{eqnarray*}
where $k \ge 0$ (for $k=0$ we adopt the standard convention $\bigcup_{j=1}^0 \ldots = \emptyset$), $S_0 \cup \ldots \cup S_k = \{1, \ldots, n \}$, 
$\pi:  \mathbb{P}^1 \to \mathbb{P}^1$ is the automorphism of $\mathbb{P}^1$ such that $\pi( p_{n+1})=0$, $\pi( p_{n+2})=1$, $\pi( p_{n+3})=\infty$,
$\mathrm{id}: \mathbb{P}^1 \to \mathbb{P}^1$ is the identity on $\mathbb{P}^1$ and $\mathrm{pt}: \mathbb{P}^1 \to \mathbb{P}^1$ collapses 
$\mathbb{P}^1$ to a point. 

By definition, $\Phi$ is injective, it is surjective onto the open part $M_{0,n}(\mathbb{P}^1,1) \subset \monp$ parameterizing stable maps with smooth 
domain $\mathbb{P}^1$ and for every subset $S \subset  \{ 1, \ldots, n \}$ with $2 \le \vert S \vert \le n$ the image 
$\Phi(\Delta_{S} \cap U_0)$ is dense in $B_S$, so that every boundary component of $\monp$  is dominated by $\Phi$. 
It follows that $\Phi$ induces a birational map $\overline{M}_{0,n+3} \dashrightarrow \monp$ which is surjective in codimension one. 

\qed


\section{}\label{ample}

According to \cite{CHS:09}, the ample cone of $\monp$ can be described in terms of natural maps: 
\begin{eqnarray*}
&\alpha:& \overline{M}_{0,n+1} \to \monp \\
&\beta_i:& \mathbb{P}^1 \to \monp, \hspace{0.2cm} i=1, \ldots, n
\end{eqnarray*}
defined in \cite{CHS:09}, 2.1 and 2.2. Indeed, by \cite{CHS:09}, Theorem 2.3, a divisor $H$ on $\monp$ is ample if and only if $\alpha^*H$ is ample on 
$\overline{M}_{0,n+1}$ and $\beta_i^*H$ is ample on $\mathbb{P}^1$ for $i=1, \ldots, n$.

In addition to the divisors $\Delta_S \subset \monbar$ and $B_S \subset \monp$, for $i=1, \ldots, n$ we introduce also the classes
$$
\psi_i := c_1(T^*_i), 
$$
where $T^*_i$ is the line bundle on $\monbar$ whose fiber over $(C, p_1, \ldots, p_n)$ is the cotangent space $(T_{p_i}C)^*$, and 
$$
L_i := \{ (C, p_1, \ldots, p_n, f: C \to \mathbb{P}^1) \in \monp: f(p_i)=0 \}. 
$$

The pullback of the classes $B_S$ and $L_i$ under the maps $\alpha$ and $\beta_i$ is computed in \cite{CHS:09}, Proposition 2.5 (see also \cite{CHS:09}, Table 1),
in terms of the classes $\Delta_S$ and $\psi_i$, namely: 
\begin{eqnarray*}
\alpha^*B_S &=& \Delta_S \hspace{0.2cm} \textrm{ if } \vert S \vert \le n-1 \\
\alpha^*B_S &=& - \psi_{n+1} \hspace{0.2cm} \textrm{ if } \vert S \vert = n \\
\alpha^*L_i &=& 0 \hspace{0.2cm} \textrm{ for every } i=1, \ldots, n \\
\beta_i^*B_S &=& \mathcal{O}_{\mathbb{P}^1}(-1) \hspace{0.2cm} \textrm{ for } S = \{1, \ldots, n \} \textrm{ and } S = \{ i \}^c \\
\beta_i^*B_S &=& 0 \hspace{0.2cm} \textrm{ otherwise } \\
\beta_i^*L_i &=& \mathcal{O}_{\mathbb{P}^1}(1) \\
\beta_j^*L_i &=& 0 \hspace{0.2cm} \textrm{ for every } j \ne i.
\end{eqnarray*}

The canonical class of $\monp$ is 
$$
K_n = -2L + \sum_{s=3}^n (s-2) B[s]
$$
where 
\begin{eqnarray*}
L &:=& \sum_{i=1}^s L_i \\
B[s] &:=& \sum_{\vert S \vert = s} B_S, \hspace{0.2cm} 2 \le s \le n
\end{eqnarray*}
(see \cite{HT:03}, p. 596). Hence we have
\begin{eqnarray*}
\alpha^*K_n &=& - (n-2) \psi_{n+1} + \sum_{s=3}^{n-1}(s-2) \sum_{ {\vert S \vert = s} \atop {n+1 \notin S}} \Delta_S\\
\beta_i^*K_n &=&  \mathcal{O}_{\mathbb{P}^1}(-2) \otimes  \mathcal{O}_{\mathbb{P}^1}(-(n-2)) \otimes \mathcal{O}_{\mathbb{P}^1}(-(n-3)) \\
&=&  \mathcal{O}_{\mathbb{P}^1}(-(2n-3)).
\end{eqnarray*}

According to Fulton's conjecture (see \cite{GKM:02}, Conjecture 0.2), a divisor on $\monbar$ is ample if and only if it has positive intersection with all one-dimensional strata, 
parameterizing $n$-pointed rational curves with at least $n-4$ singular points. More explicitly, let
$$
H = \sum_{ \vert S \vert \ge 1} c_S \Delta_S,
$$
where we adopt the convention $\Delta_{\{i\}} := - \psi_i$ for every $i=1, \ldots, n$. By \cite{GKM:02}, Theorem 2.1, the divisor $H$ has positive intersection 
with all one-dimensional strata if and only if  
$$
c_{I \cup J} + c_{I \cup K} + c_{I \cup L} - c_I - c_J - c_K - c_L > 0
$$
for every partition $I \cup J \cup K \cup L$ of $\{1, \ldots, n \}$. 

By \cite{KM:96}, Theorem 1.2(3), Fulton's conjecture holds for $n \le 7$.


\section{}\label{fano}

Finally we are going to check that $\monp$ is a log Fano variety (hence a Mori dream space) for $n \le 5$ but not for $n=6$.

\begin{Lemma}\label{n=4}
On $\overline{M}_{0,4}(\mathbb{P}^1,1)$ the divisor $K_4 + B[4]$ is anti-ample, hence $\overline{M}_{0,4}(\mathbb{P}^1,1)$ is log Fano.
\end{Lemma}

\proof By \cite{CHS:09}, Proposition 2.5, we have 
\begin{eqnarray*}
\alpha^*(K_4 + B[4]) &=& - 3 \psi_5 + \sum_{ {\vert S \vert = 3} \atop {5 \notin S}} \Delta_S\\
\beta_i^*(K_4 + B[4]) &=&  \mathcal{O}_{\mathbb{P}^1}(-2) \otimes  \mathcal{O}_{\mathbb{P}^1}(-3) \otimes \mathcal{O}_{\mathbb{P}^1}(-1) \\
&=&  \mathcal{O}_{\mathbb{P}^1}(-6).
\end{eqnarray*}
It is clear that $\beta_i^*(K_4 + B[4])$ is anti-ample on $\mathbb{P}^1$; on the other hand, in order to check that $\alpha^*(K_4 + B[4])$ is anti-ample 
on $\overline{M}_{0,5}$, by \cite{GKM:02}, Theorem 2.1, we have to consider the following partitions $I \cup J \cup K \cup L$ of $\{1, \ldots, 5 \} = 
\{ a,b,c,d,5 \}$:
\begin{itemize}
\item $\{ a \} \cup \{ b \} \cup \{ c \} \cup \{ d,5 \}$ 
\item $\{ a \} \cup \{ b \} \cup \{ 5 \} \cup \{ c,d \}$.
\end{itemize}
If
$$
\alpha^*(K_4 + B[4]) = \sum_{\vert S \vert \ge 1} c_S \Delta_S
$$
then
$$
c_{I \cup J} + c_{I \cup K} + c_{I \cup L} - c_I - c_J - c_K - c_L = -1
$$
in both cases listed above, hence $\alpha^*(K_4 + B[4])$ is anti-ample by \cite{KM:96}, Theorem 1.2(3), and $K_4 + B[4]$ is anti-ample by \cite{CHS:09},Theorem 2.3.

\qed

\begin{Lemma}\label{n=5}
Let $D = \frac{1}{4} B[2]+ \frac{1}{4} B[4]+B[5]$ on $\overline{M}_{0,5}(\mathbb{P}^1,1)$. The divisor $K_5 + D$ is anti-ample, 
hence $\overline{M}_{0,5}(\mathbb{P}^1,1)$ is log Fano.
\end{Lemma}

\proof By \cite{CHS:09}, Proposition 2.5, we have 
\begin{eqnarray*}
\alpha^*(K_5 + D) &=& - 4 \psi_6 + \frac{1}{4} \sum_{ {\vert S \vert = 2} \atop {6 \notin S}} \Delta_S  
+ \sum_{ {\vert S \vert = 3} \atop {6 \notin S}} \Delta_S + \left(  2 + \frac{1}{4} \right) \sum_{ {\vert S \vert = 4} \atop {6 \notin S}} \Delta_S  \\
\beta_i^*(K_5 + D) &=&  \mathcal{O}_{\mathbb{P}^1}(-2) \otimes  \mathcal{O}_{\mathbb{P}^1}(-4) \otimes \mathcal{O}_{\mathbb{P}^1} \left( -2- \frac{1}{4} \right) \\
&=&  \mathcal{O}_{\mathbb{P}^1} \left(-8 - \frac{1}{4} \right).
\end{eqnarray*}
It is clear that $\beta_i^*(K_5 + D)$ is anti-ample on $\mathbb{P}^1$; on the other hand, in order to check that $\alpha^*(K_5 + D)$ is anti-ample 
on $\overline{M}_{0,6}$, by \cite{GKM:02}, Theorem 2.1, we have to consider the following partitions $I \cup J \cup K \cup L$ of $\{1, \ldots, 6 \} = 
\{ a,b,c,d,e,6 \}$:
\begin{itemize}
\item $\{ a \} \cup \{ b \} \cup \{ c \} \cup \{ d,e,6 \}$
\item $\{ a \} \cup \{ b \} \cup \{ 6 \} \cup \{ c,d,e \}$
\item $\{ a \} \cup \{ b \} \cup \{ c,d \} \cup \{ e,6 \}$
\item $\{ a \} \cup \{ 6 \} \cup \{ b,c \} \cup \{ d,e \}$.
\end{itemize}
If
$$
\alpha^*(K_5 + D) = \sum_{\vert S \vert \ge 1} c_S \Delta_S
$$
then 
$$
c_{I \cup J} + c_{I \cup K} + c_{I \cup L} - c_I - c_J - c_K - c_L = -\frac{1}{4}
$$
in all cases listed above, hence $\alpha^*(K_5 + D)$ is anti-ample by \cite{KM:96}, Theorem 1.2(3), and $K_5 + D$ is anti-ample by \cite{CHS:09},Theorem 2.3.

\qed

\begin{Lemma}\label{n=6}
Let $D = a_2 B[2]+ a_3 B[3]+ a_4 B[4]+ a_5 B[5] + a_6 B[6]$ on $\overline{M}_{0,6}(\mathbb{P}^1,1)$ with $a_i \in \mathbb{Q}$. 
If $a_4 \ge 0$ and $a_6 \le 1$ then $K_6 + D$ is not anti-ample.
\end{Lemma}

\proof By \cite{CHS:09}, Proposition 2.5, we have 
$$
\alpha^*(K_6 + D) = \sum_{\vert S \vert \ge 1} c_S \Delta_S
$$
with 
\[
c_S =
\begin{cases}
4 + a_6 \hspace{0.2cm} \textrm{ if } S = \{ 7 \} \\
0 \hspace{0.2cm} \textrm{ if } \vert S \vert = 1,  S \ne \{ 7 \} \\ 
3 + a_5 \hspace{0.2cm} \textrm{ if } \vert S \vert = 2,  7 \in S \\
a_2 \hspace{0.2cm} \textrm{ if } \vert S \vert = 2,  7 \notin S \\
2 + a_4 \hspace{0.2cm} \textrm{ if } \vert S \vert = 3,  7 \in S \\
1 + a_3 \hspace{0.2cm} \textrm{ if } \vert S \vert = 3,  7 \notin S. \\
\end{cases}
\]
Consider the following partitions of $\{1, \ldots, 7 \} = \{ a,b,c,d,e,f,7 \}$:

(i) $\{ a \} \cup \{ b \} \cup \{ c \} \cup \{ d,e,f,7 \}$

(ii) $\{ a \} \cup \{ b \} \cup \{ c,d \} \cup \{ e,f,7 \}$

(iii) $\{ 7 \} \cup \{ a,b \} \cup \{ c,d \} \cup \{ e,f \}$.

\noindent According to \cite{GKM:02}, Theorem 2.1, the corresponding necessary conditions for $\alpha^*(K_6 + D)$ to be anti-ample are:

(i) $3 a_2 - a_3 - 1 < 0$ 

(ii) $2 a_3 - a_4 < 0$

(iii) $3 a_4 - 3 a_2 - a_6 + 2 < 0.$  

\noindent
Hence we deduce: 

(ii) $a_3 < \frac{a_4}{2}$

(i) $a_2 < \frac{1}{3}+\frac{a_3}{3}< \frac{1}{3}+\frac{a_4}{6}$

(iii) $a_2 > \frac{2-a_6}{3} + a_4$

\noindent
which is impossible if  $a_4 \ge 0$ and $a_6 \le 1$.

\qed

\begin{Corollary}\label{dream}
If $n \le 5$ then $\monp$ is a Mori dream space. 
\end{Corollary}

\proof
If $n \le 3$ we exploit the isomorphism $\monp \cong \mathbb{P}^1[n]$, where $\mathbb{P}^1[n]$ denotes the Fulton-MacPherson compactification 
(see \cite{Ma:17}, Remark 1.4)  and the fact that $\mathbb{P}^1[n]$ is Fano for $n \le 3$  (see \cite{Ma:17}, p. 108). If $n=4,5$ then 
$\monp$ is log Fano by Lemma \ref{n=4} and Lemma \ref{n=5}. Hence $\monp$ is a Mori dream space 
for $n \le 5$ by \cite{BCHM:10}, Corollary 1.3.2. 

\qed

\begin{Remark}\label{nohope} By Lemma \ref{n=6}, there is no hope to deduce from \cite{BCHM:10}, Corollary 1.3.2, that 
$\overline{M}_{0,6}(\mathbb{P}^1,1)$ is a Mori dream space. 
\end{Remark}

\end{document}